\documentclass{amsart}
\usepackage{setspace}
\usepackage{a4}
\usepackage{amssymb,amsmath,amsthm,latexsym}
\usepackage{amsfonts}
\usepackage{graphicx}
\usepackage{textcomp}
\usepackage{cite}
\usepackage{enumerate}
\usepackage[mathscr]{euscript}
\usepackage{mathtools}
\newtheorem{theorem}{Theorem}[section]

\newtheorem{corollary}[theorem] {Corollary}
\newtheorem{definition}[theorem]{Definition}
\newtheorem{example}[theorem]{Example}
\newtheorem{lemma} [theorem]{Lemma}

\newtheorem{problem}[theorem]{Problem}
\newtheorem{proposition}[theorem]{Proposition}

\newtheorem{motivation}[theorem] {Motivation}
\title{This is the title}
\usepackage{hyperref}
\usepackage{enumerate}
\usepackage{mathtools}
\usepackage{amsmath}
\usepackage{tikz}
\usepackage{graphicx}
\usepackage{textcomp}
\usetikzlibrary{cd}
\usepackage{fancyhdr}

\begin{document}
\begin{center}
	\hrule 	\hrule 	\hrule 	\hrule \hrule 
	\vspace{0.3cm}
{\bf{LIPSCHITZ P-APPROXIMATE SCHAUDER FRAMES}}\\
\vspace{0.3cm}
\hrule\hrule 	\hrule 	\hrule \hrule 
\vspace{0.3cm}
{\bf{K. MAHESH KRISHNA}} \\
Statistics and Mathematics Unit\\
Indian Statistical Institute, Bangalore Centre\\
Karnataka 560 059 India\\
Email: kmaheshak@gmail.com \\
and \\
{\bf{P. SAM JOHNSON}} \\
Department of Mathematical and Computational Sciences \\
National Institute of Technology Karnataka, Surathkal \\
Mangaluru 575 025, India \\
Email: sam@nitk.edu.in

Date: \today
\end{center}

\hrule
\vspace{0.3cm}
%--------------------------------------
\textbf{Abstract}: With the aim of representing subsets of Banach spaces as an infinite  series using Lipschitz functions, we study a variant of metric frames which we call  Lipschitz p-approximate Schauder frames (Lipschitz p-ASFs). We characterize Lipschitz p-ASFs  and their duals completely using the canonical Schauder basis for classical sequence spaces. Similarity of  Lipschitz p-ASF is   introduced and characterized.

\textbf{Keywords}:  Frame, Lipschitz map, bi-Lipschitz map.

\textbf{Mathematics Subject Classification (2020)}: 42C15, 26A16.
\vspace{0.3cm}
\hrule
\tableofcontents
\hrule
\section{Introduction}

Grochenig  in 1991 introduced the notion of Banach frames \cite{GROCHENIG} as a generalization of notion of frames for Hilbert spaces introduced by  Duffin and Schaeffer in  1952 \cite{DUFFIN}. This notion originated from the study of atomic decompositions and coorbit spaces arising from square integrable representations of locally compact groups developed by  Feichtinger and Grochenig in 1980's   \cite{FEICHTINGERGROCHENIG, FEICHTINGERGROCHENIG2, FEICHTINGERGROCHENIG1}. Casazza, Han and Larson in 2000 explored the connection between Banach frames and atomic decompositions and  introduced the notion of (unconditional) Schauder frames \cite{CASAZZAHANLARSON}.  In 2001, Aldroubi, Sun and Tang introduced the notion of p-frames and p-Riesz bases for Banach spaces, $1\leq p <\infty$ \cite{ALDROUBISUNTANG}. These notions have  been  generalized by Casazza, Christensen and  Stoeva by introducing the notion of $\mathcal{X}_d$-frames \cite{CHRISTENSENSTOEVA, CASAZZACHRISTENSENSTOEVA}. A slight variant notions of  $\mathcal{X}_d$-frames  for Banach spaces was given by Terekhin  \cite{TEREKHIN1, TEREKHIN2, TEREKHIN3}.  In 2014, Thomas, Freeman, Odell,  Schlumprecht, and Zsak  \cite{FREEMANODELL, THOMAS} introduced the notion of approximate Schauder  frames as a generalization of notion of Schauder frames by Casazza, Dilworth, Odell, Schlumprecht and Zsak \cite{CASAZZADILWORTHODELL} (also see \cite{DAUBECHIESDEVORE}). In 2021, Krishna and Johnson characterized some classes of approximate Schauder frames \cite{KRISHNAJOHNSON}. In 2021, Krishna introduced  Lipschitz atomic decompositions for subsets of Banach spaces \cite{KRISHNA4}.  In 2022,  Krishna and Johnson introduced metric frames which has surprising connections with subsets of Banach spaces using Lipschitz-free Banach spaces \cite{KRISHNAJOHNSON2}. We now ask the following question which is the main motivation for writing the paper.
\begin{motivation}\label{MOT}
\textbf{Can we represent a subset (which need not be a subspace) of a Banach space as an infinite series using Lipschitz maps and elements of the set?}
\end{motivation}
Note that we can not demand linear functionals in the above problem as we are not considering subspaces. Motivated from \ref{MOT} we study representation of subsets (need not be subspaces) of Banach spaces using Lipschitz functions.

The paper is organized  as follows. We introduce the notion fo Lipschitz p-approximate Schauder frame (Lipschitz p-ASF) for subsets of Banach spaces in Definition \ref{LIPASFDEFINITION}. Followed by  interesting Examples \ref{FEXAMPLE}, \ref{SEXAMPLE} and \ref{PRO}, factorization property of Lipschitz frame map is derived in Theorem \ref{FP}. Lipschitz p-ASFs are characterized in Theorem \ref{LASFCHAR}. Next we introduce the notion of dual frames in Definition \ref{DUALDEF} and classify them in Theorem \ref{ALLDUAL}. Definition \ref{SIMDEF} introduces the notion of similarity and Theorem \ref{SEQUENTIALSIMILARITY}  gives an operator-theoretic characterization for similarity. Orthogonality of frames is introduced in Definition \ref{ORTHOGONALDEF} and interpolation result is derived in Theorem \ref{INTERPOLATION}. We end by formulating  open Problem \ref{PROBLEM}.

 \section{Lipschitz p-approximate Schauder frames}\label{SECTION2}
  Let $\mathcal{X}$ be  a real or complex Banach space be a     $\mathcal{M}$ is a non-empty subset of $\mathcal{X}$.  The identity operator on $\mathcal{M}$  is denoted by $I_\mathcal{M}$.  The set of all Lipschitz functions from $\mathcal{M}$ to $\mathcal{X}$ is denoted by $\operatorname{Lip}(\mathcal{M}, \mathcal{X})$. For $1\leq p<\infty$, the canonical Schauder basis for $\ell^p(\mathbb{N})$ is  denoted by $\{e_n\}_n$ and its coordinate functionals are denoted by $\{\zeta_n\}_n$.  We introduce the following important notion as a first step in answering Motivation \ref{MOT}.
  \begin{definition}\label{LIPASFDEFINITION}
 	For $1\leq p<\infty$, let $\mathcal{X}$ be a  Banach space and 	$\mathcal{M}$ be a subset (need not be a subspace) of $\mathcal{X}$. Let  $\{\tau_n\}_n$ be a sequence in  $\mathcal{M}$ and 	$\{f_n\}_n$ be a sequence in  $\operatorname{Lip}(\mathcal{M}, \mathcal{X})$.  The pair  $ (\{f_n \}_{n}, \{\tau_n \}_{n}) $  	is said to be a \textbf{Lipschitz p-approximate Schauder frame} (we write Lipschitz p-ASF)  if the following conditions hold. 
 	\begin{enumerate}[\upshape(i)]
 		\item The map (\textbf{analysis map}) $ \theta_f: \mathcal{M}\ni x \mapsto \theta_f x\coloneqq \{f_n(x)\}_n \in \ell^p(\mathbb{N})$ is a  well-defined  Lipschitz map.
 	 		\item The map (\textbf{synthesis operator}) $
 	 		\theta_\tau : \ell^p(\mathbb{N}) \ni \{a_n\}_n \mapsto \theta_\tau \{a_n\}_n\coloneqq \sum_{n=1}^\infty a_n\tau_n \in \mathcal{X}	$ is a  well-defined bounded linear operator. 
 	 			\item The map (\textbf{Lipschitz frame  map}) $S_{f, \tau}:\mathcal{M}\ni x \mapsto S_{f, \tau}x\coloneqq \sum_{n=1}^\infty f_n(x)\tau_n \in\mathcal{M}$ is a well-defined invertible bi-Lipschitz map and 
 	 				\begin{align}\label{REPBANACH0}
 	 				x=\sum_{n=1}^\infty f_n(x)S_{f, \tau}^{-1}\tau_n, \quad \forall x \in
 	 				\mathcal{M}.
 	 			\end{align}
 	\end{enumerate}
  If $S_{f, \tau}=I_\mathcal{M}$, then we say that $ (\{f_n \}_{n}, \{\tau_n \}_{n}) $ is a \textbf{Lipschitz p-Schauder frame} (we write Lipschitz p-SF).     If we do not impose the condition `invertible bi-Lipschitz' and Equation (\ref{REPBANACH0}) in  (iii), then we say that $ (\{f_n \}_{n}, \{\tau_n \}_{n}) $  	is  a \textbf{Lipschitz p-Bessel sequence} (we write Lipschitz p-BS) for $\mathcal{M}$. 
  \end{definition}
Whenever $\mathcal{M}=\mathcal{X}$,  and $f_n$'s are all linear, Definition \ref{LIPASFDEFINITION}  reduces to definition of p-ASF given in \cite{KRISHNAJOHNSON}. It is important to note that  the partial sums of series  in (iii)  of  Definition \ref{LIPASFDEFINITION}  need not  be inside $\mathcal{M}$ (which may not be as it is only a subset) but only demanding limit has to be inside $\mathcal{M}$.  Definition \ref{LIPASFDEFINITION} says that there are $a,b,c,d>0$ satisfying following:
  \begin{align*}
 	&a\|x-y\|\leq \left\|\sum_{n=1}^\infty
 	(f_n(x)-f_n(y))\tau_n \right\|\leq b\|x-y\|, \quad \forall x,y \in  \mathcal{M},\\
 	& 	\left(\sum_{n=1}^{\infty}|f_n(x)-f_n(y)|^p	\right)^\frac{1}{p}\leq c\|x-y\|, \quad \forall x, y \in \mathcal{M}, \\
 	& 	\left\|\sum_{n=1}^\infty a_n\tau_n\right\|\leq d \left(\sum_{n=1}^{\infty}|a_n|^p	\right)^\frac{1}{p}, \quad \forall \{a_n\}_n \in \ell^p(\mathbb{N}).
 \end{align*}  
 We call $a$ as lower Lipschitz frame bound, $b$ as upper Lipschitz frame bound, $c$ as Lipschitz analysis bound and $d$ as Lipschitz synthesis bound.
 We give various interesting examples of Lipschitz p-ASFs.
 \begin{example}\label{FEXAMPLE}
 	Let $\mathcal{X}\coloneqq \mathbb{C}$, $p=1$ and 
 	\begin{align*}
 		\mathcal{M}\coloneqq \left\{z \in \mathbb{C}: |z|\leq \frac{1}{2}|z+1|\right\}
 		=\left\{x+iy:x, y  \in \mathbb{R},  \left(x-\frac{1}{3}\right)^2+y^2\leq \left(\frac{2}{3}\right)^2\right\}.
 	\end{align*}
 For $n\in \mathbb{N}$, define 
 \begin{align*}
 	f_n: \mathcal{M} \ni z \mapsto f_n(z)\coloneqq \left(\frac{z}{1+z}\right)^n \in \mathbb{C}, \quad \tau_n\coloneqq 1 \in \mathcal{M}.
 \end{align*}
We first show that $f_n$ is Lipschitz for all $n$. For $z \in \mathcal{M}$, 
\begin{align*}
	1-|z+1|\leq \big|1-|z+1|\big|\leq \big|1-(z+1)\big|=|z|\leq \frac{1}{2}|z+1|.
\end{align*}
Hence 
\begin{align*}
	|z+1|\geq \frac{2}{3}, \quad \forall z \in \mathcal{M}.
\end{align*}
Let $z, w \in \mathcal{M}$.  Then for each $n \in \mathbb{N}$,
\begin{align*}
&|f_n(z)-f_n(w)|=\left|\left(\frac{z}{1+z}\right)^n-\left(\frac{w}{1+w}\right)^n \right|\\
&=\left|\frac{z}{1+z}-\frac{w}{1+w}\right|\left|\left(\frac{z}{1+z}\right)^{n-1}+\left(\frac{z}{1+z}\right)^{n-2}\left(\frac{w}{1+w}\right)+\cdots +\left(\frac{z}{1+z}\right)\left(\frac{w}{1+w}\right)^{n-2}+\left(\frac{w}{1+w}\right)^{n-1}\right|\\
&\leq \frac{|z-w|}{|1+z||1+w|}\frac{n}{2^{n-1}}\leq \frac{9}{4}\frac{n}{2^{n-1}}|z-w|.
\end{align*}
Therefore each $f_n$ is Lipschitz. Set 
\begin{align*}
	r\coloneqq \sum_{n=1}^{\infty} \frac{n}{2^{n-1}}<\infty.
\end{align*}
We then  see that for $z, w \in \mathcal{M}$, 
\begin{align*}
	\|\theta_f z-\theta_f w\|&=\sum_{n=1}^{\infty}|f_n(z)-f_n(w)|=\sum_{n=1}^{\infty}\left|\left(\frac{z}{1+z}\right)^n-\left(\frac{w}{1+w}\right)^n \right|\\
	&\leq \sum_{n=1}^{\infty}\frac{9}{4}\frac{n}{2^{n-1}}|z-w|=\frac{9}{4}r|z-w|.
\end{align*}
Therefore $\theta_f$ is Lipschitz. Clearly 
\begin{align*}
	\theta_\tau:\ell^1(\mathbb{N}) \ni \{a_n\}_n \mapsto \sum_{n=1}^{\infty}a_n\cdot 1 \in \mathbb{C}
\end{align*}
is a well-defined bounded linear operator. Finally,  we observe that for $z \in \mathcal{M}$, we have $\frac{|z|}{|z+1|}<1$  and hence 
\begin{align*}
S_{f, \tau}	z=\sum_{n=1}^{\infty} f_n(z)\tau_n=\sum_{n=1}^{\infty}\left(\frac{z}{1+z}\right)^n\cdot 1=\frac{1}{1-\frac{z}{1+z}}-1=z, \quad \forall z \in \mathcal{M}.
\end{align*}
Thus we proved that $ (\{f_n\}_{n}, \{\tau_n\}_{n}) $  is a Lipschitz 1-SF for $\mathcal{M}$. 
 \end{example}
 \begin{example}\label{SEXAMPLE}
	Let $\mathcal{X}\coloneqq \mathbb{R}$, $p=1$ and   $\mathcal{M}\coloneqq[1,\infty)$ 
	 For $n\in \mathbb{N}\cup \{0\}$, define  $f_n:\mathcal{M}\to \mathbb{R}$ by 
\begin{align*}
	f_0(x)&\coloneqq 1, \quad \forall x \in \mathcal{M},\\
	f_n(x)&\coloneqq \frac{(\log x)^n}{n!}, \quad \forall x \in \mathcal{M}, \forall n\geq 1 
\end{align*}
and $\tau_n \coloneqq 1 \in \mathcal{M}$. 
Then $f_n'(x)=\frac{(\log x)^{(n-1)}}{(n-1)!x}$, $\forall x \in \mathcal{M}, \forall n\geq1.$ Since   $f_n'$ is bounded on $\mathcal{M}$ for all  $ n\geq1$,  $f_n$ is Lipschitz on $\mathcal{M}$ for all  $n\geq1$. For $x, y \in \mathcal{M},$ with $x<y$, we see that 
\begin{align*}
		\|\theta_f x-\theta_f y\|=\sum_{n=0}^{\infty}|f_n(x)-f_n(y)|=\sum_{n=0}^{\infty}\frac{(\log
		y)^n}{n!}-\sum_{n=0}^{\infty}\frac{(\log x)^n}{n!}
	=e^{\log y}-e^{\log x}=y-x=|x-y|.
\end{align*}
Therefore $\theta_f$ is Lipschitz. It is clear that  $\theta_\tau$  is a well-defined bounded linear operator. For $x \in \mathcal{M}$, 
\begin{align*}
	S_{f, \tau}	x=\sum_{n=1}^{\infty} f_n(x)\tau_n=\sum_{n=0}^{\infty}\frac{(\log x)^n}{n!}\cdot 1=x.
\end{align*}
Hence $ (\{f_n\}_{n}, \{\tau_n\}_{n}) $ is a Lipschitz  1-SF for $\mathcal{M}$.
 \end{example}
 \begin{example}\label{PRO}
  	For $1\leq p<\infty$, let $\mathcal{X}$ be a  Banach space and $\mathcal{M}$ be a subset of $\mathcal{X}$.   Assume that there is a Lipschitz map $U:\mathcal{M} \rightarrow \ell^p(\mathbb{N})$, a bounded linear operator $ V:\ell^p(\mathbb{N})\to \mathcal{X}$  such that $VU(\mathcal{M})\subseteq \mathcal{M}$, $Ve_n\in \mathcal{M} $ for all $n \in \mathbb{N}$,   $VU:\mathcal{M} \to \mathcal{M}$ is an  invertible bi-Lipschitz  map and  
  	\begin{align*}
  		x=\sum_{n=1}^{\infty}\zeta_n(Ux)(VU)^{-1}Ve_n, \quad \forall x \in \mathcal{M}.
  	\end{align*}
  	Let $\{e_n\}_n$ denote the canonical Schauder basis for  $\ell^p(\mathbb{N})$  and let $\{\zeta_n\}_n$ denote the coordinate functionals associated with $\{e_n\}_n$.	Define 
  	\begin{align*}
  		f_n\coloneqq \zeta_n U, \quad \tau_n\coloneqq Ve_n, \quad \forall n \in \mathbb{N}.
  	\end{align*}
  	Then $ (\{f_n\}_{n}, \{\tau_n\}_{n}) $ is Lipschitz p-ASF for 	$\mathcal{M}$. If $VU=I_\mathcal{M}$, then  $(\{f_n\}_{n}, \{\tau_n\}_{n}) $ is a Lipschitz p-SF for $\mathcal{M}$. 
  \end{example}
 We show in the sequel that (in Theorem \ref{LASFCHAR}) every Lipschitz p-ASF can be written in the form of Example \ref{PRO}. Following theorem gives various fundamental factorization properties of Lipschitz p-ASFs whose proof is a direct calculation.
  \begin{theorem}\label{FP}
  	Let $ (\{f_n \}_{n}, \{\tau_n \}_{n}) $ be a Lipschitz p-ASF for $\mathcal{M}\subseteq \mathcal{X}$.  Then
  	\begin{enumerate}[\upshape(i)]
  		\item We have 
  		\begin{align}\label{REPBANACH}
  			x=\sum_{n=1}^\infty (f_nS_{f, \tau}^{-1})(x) \tau_n, \quad \forall x \in
  			\mathcal{M}.
  		\end{align}
  		\item $ (\{f_nS_{f, \tau}^{-1} \}_{n}, \{S_{f, \tau}^{-1} \tau_n \}_{n}) $ is a  Lipschitz p-ASF for $\mathcal{M}$.
  		\item The analysis map $\theta_f$	is injective. 
  		\item 
  		The synthesis operator  $\theta_\tau$	is surjective.
  		\item Lipschitz frame map $S_{f, \tau}$ factors  as $S_{f, \tau}=\theta_\tau\theta_f.$
  		\item  $P_{f, \tau}\coloneqq\theta_fS_{f,\tau}^{-1}\theta_\tau:\ell^p(\mathbb{N})\to \ell^p(\mathbb{N})$ is a Lipschitz projection onto   $\theta_f(\mathcal{M})$.
  	\end{enumerate}
  \end{theorem}
  Holub  characterized  frames for Hilbert spaces using standard orthonormal basis for the standard Hilbert space \cite{HOLUB}. This result has been derived for Banach spaces in \cite{KRISHNAJOHNSON}. We show that such a result can be derived for Lipschitz p-ASFs.  
  \begin{theorem}\label{LASFCHAR}
  	A pair  $ (\{f_n\}_{n}, \{\tau_n\}_{n}) $ is a Lipschitz p-ASF for 	$\mathcal{M}\subseteq \mathcal{X}$
  	if and only if 
  	\begin{align*}
  		f_n=\zeta_n U, \quad \tau_n=Ve_n, \quad \forall n \in \mathbb{N},
  	\end{align*}  
  where $U:\mathcal{M} \rightarrow \ell^p(\mathbb{N})$ is a Lipschitz map, $ V: \ell^p(\mathbb{N})\to \mathcal{X}$ is a  bounded linear operator such that $VU(\mathcal{M})\subseteq \mathcal{M}$, $Ve_n\in \mathcal{M} $ for all $n \in \mathbb{N}$,  $VU:\mathcal{M} \to \mathcal{M}$ is an invertible bi-Lipschitz map and 
  \begin{align*}
  	x=\sum_{n=1}^{\infty}\zeta_n(Ux)(VU)^{-1}Ve_n, \quad \forall x \in \mathcal{M}.
  \end{align*}
  \end{theorem}
  \begin{proof}
  	$(\Leftarrow)$ Clearly $\theta_f$ is Lipschitz and $\theta_\tau$ is a  bounded linear operator. Now let $x\in \mathcal{M}$. Then 
  	\begin{align}\label{PASFFIRSTTHEOREMEQUATION}
  		S_{f, \tau}x= \sum_{n=1}^\infty
  		f_n(x)\tau_n=\sum_{n=1}^\infty \zeta_n(Ux)Ve_n=V\left(\sum_{n=1}^\infty \zeta_n(Ux)e_n\right)=VUx.
  	\end{align} 
  	Hence $S_{f, \tau}$ is an  invertible bi-Lipschitz map. \\
  	$(\Rightarrow)$ Define $U\coloneqq \theta_f$, $V\coloneqq \theta_\tau$. Then $(\zeta_nU)(x)=(\zeta_n\theta_f)(x)=\zeta_n(\{f_k(x)\}_k)=f_n(x)$, $\forall x \in \mathcal{M}$, $Ve_n=\theta_\tau e_n=\tau_n$, $\forall n \in \mathbb{N}$ and $VU=\theta_\tau \theta_f=S_{f, \tau}$ which is an invertible bi-Lipschitz map.
  \end{proof}
\begin{corollary}
	\begin{enumerate}[\upshape(i)]
	\item 	A pair  $ (\{f_n\}_{n}, \{\tau_n\}_{n}) $ is a Lipschitz p-SF for 	$\mathcal{M}\subseteq \mathcal{X}$  if and only if $	f_n=\zeta_n U,  \tau_n=Ve_n,  \forall n \in \mathbb{N}, $ where $U:\mathcal{M} \rightarrow \ell^p(\mathbb{N})$ is a Lipschitz map, $ V:  \ell^p(\mathbb{N})\to \mathcal{X}$ is a  bounded linear operator such that $VU(\mathcal{M})\subseteq \mathcal{M}$, $Ve_n\in \mathcal{M} $ for all $n \in \mathbb{N}$ and $VU=I_\mathcal{M}$.
	\item  A pair  $ (\{f_n\}_{n}, \{\tau_n\}_{n}) $ is a Lipschitz p-BS for 	$\mathcal{M}\subseteq \mathcal{X}$  if and only if $	f_n=\zeta_n U,  \tau_n=Ve_n,  \forall n \in \mathbb{N}, $ where $U:\mathcal{M} \rightarrow\ell^p(\mathbb{N})$ is a Lipschitz map, $ V: \ell^p(\mathbb{N})\to \mathcal{X}$ is a  bounded linear operator such that $VU(\mathcal{M})\subseteq \mathcal{M}$ and  $Ve_n\in \mathcal{M} $ for all $n \in \mathbb{N}$. 
	\end{enumerate}	
\end{corollary}
Equations (\ref{REPBANACH0}) and  (\ref{REPBANACH}) lead  us to define the notion of dual frame as follows.
\begin{definition}\label{DUALDEF}
	Let $ (\{f_n\}_{n}, \{\tau_n\}_{n}) $ be a Lipschitz p-ASF for 	$\mathcal{M}\subseteq\mathcal{X}$. 	A Lipschitz p-ASF $ (\{g_n \}_{n}, \{\omega_n \}_{n}) $ for $\mathcal{M}\subseteq\mathcal{X}$ is said to be  a \textbf{dual}   for $ (\{f_n \}_{n}, \{\tau_n \}_{n}) $ if 
	\begin{align*}
		x=\sum_{n=1}^\infty g_n(x) \tau_n=\sum_{n=1}^\infty
		f_n(x) \omega_n, \quad \forall x \in
		\mathcal{M}.
	\end{align*}
\end{definition}
We can give a  characterization  of dual frames  by using analysis map and synthesis operator.
\begin{proposition}\label{ORTHOGONALPRO}
	Given two  Lipschitz p-ASFs $ (\{f_n\}_{n}, \{\tau_n\}_{n}) $ and $ (\{g_n \}_{n}, \{\omega_n \}_{n}) $ for $\mathcal{M}\subseteq\mathcal{X}$, the following are equivalent:
	\begin{enumerate}[\upshape(i)]
		\item  $ (\{g_n \}_{n}, \{\omega_n \}_{n}) $ is a dual  for $ (\{f_n \}_{n}, \{\tau_n \}_{n}) $.
		\item $\theta_\tau\theta_g =\theta_\omega\theta_f =I_\mathcal{M}$.
	\end{enumerate}
\end{proposition}
Equations (\ref{REPBANACH0}) and (\ref{REPBANACH}) show that the Lipschitz p-ASF  $ (\{f_nS_{f, \tau}^{-1} \}_{n}, \{S_{f, \tau}^{-1} \tau_n \}_{n}) $ is a  dual for $ (\{f_n\}_{n}, \{\tau_n\}_{n}) $.  We call   $ (\{f_nS_{f, \tau}^{-1} \}_{n}, $ $ \{S_{f, \tau}^{-1} \tau_n \}_{n}) $ as the  \textbf{canonical dual} for  $ (\{f_n\}_{n}, \{\tau_n\}_{n}) $. With this notion, the following theorem  is evident.
\begin{theorem}
	Let $ (\{f_n\}_{n}, \{\tau_n\}_{n}) $ be a  Lipschitz p-ASF for $ \mathcal{M}\subseteq\mathcal{X}$ with frame bounds $ a$ and $ b.$ Then the following statements hold good:
	\begin{enumerate}[\upshape(i)]
		\item The canonical dual  for the canonical dual   for  $ (\{f_n\}_{n}, \{\tau_n\}_{n}) $ is itself.
		\item$ \frac{1}{b}, \frac{1}{a}$ are frame bounds for the canonical dual for $ (\{f_n\}_{n}, \{\tau_n\}_{n}) $.
		\item If $ a, b $ are optimal frame bounds for $ (\{f_n\}_{n}, \{\tau_n\}_{n}) $, then $ \frac{1}{b}, \frac{1}{a}$ are optimal  frame bounds for its canonical dual.
	\end{enumerate} 
\end{theorem}
In 1995 Li derived a    characterization of dual frames using standard orthonormal basis for  $\ell^2(\mathbb{N})$ \cite{LI}.   For Banach spaces, such a characterization using canonical Schauder basis for  $\ell^p(\mathbb{N})$    is derived in \cite{KRISHNAJOHNSON}. Now we derive such characterization for Lipschitz p-ASF.
\begin{lemma}\label{ASFLEMMA1}
	Let  $ (\{f_n \}_{n}, \{\tau_n \}_{n}) $  be a  Lipschitz p-ASF for   $\mathcal{M}\subseteq \mathcal{X}$. Then a  Lipschitz p-ASF  $ (\{g_n \}_{n}, \{\omega_n \}_{n}) $ for $\mathcal{M}$  is a dual  for $ (\{f_n \}_{n}, \{\tau_n \}_{n}) $ if and only if
	\begin{align*}
		g_n=\zeta_n U, \quad \omega_n=Ve_n, \quad \forall n \in \mathbb{N},
	\end{align*} 
	where $ U:\mathcal{M} \rightarrow \ell^p(\mathbb{N})$ is  a Lipschitz  right-inverse of $ \theta_\tau$ and  $V: \ell^p(\mathbb{N}) \rightarrow \mathcal{X}$ is a linear bounded left-inverse of $ \theta_f$ such that $VU(\mathcal{M})\subseteq \mathcal{M}$, $Ve_n\in \mathcal{M} $ for all $n \in \mathbb{N}$, $ VU$ is  an invertible bi-Lipschitz map and 
\begin{align*}
	x=\sum_{n=1}^{\infty}\zeta_n(Ux)(VU)^{-1}Ve_n, \quad \forall x \in \mathcal{M}.
\end{align*}
\end{lemma}
\begin{proof}
	$(\Leftarrow)$ Using  the `if' part of proof of Theorem \ref{LASFCHAR}, we get that $ (\{g_n \}_{n}, \{\omega_n \}_{n}) $ is a Lipschitz p-ASF for $\mathcal{M}$. We  check for duality of $ (\{g_n \}_{n}, \{\omega_n \}_{n}) $:  $\theta_\tau\theta_g=\theta_\tau U=I_\mathcal{M} $, $ \theta_\omega\theta_f=V\theta_f =I_\mathcal{M}$.\\
	$(\Rightarrow)$ Let $ (\{g_n \}_{n}, \{\omega_n \}_{n}) $ be a dual Lipschitz p-ASF  for  $ (\{f_n \}_{n}, \{\tau_n \}_{n}) $.  Then $\theta_\tau\theta_g =I_\mathcal{M} $, $ \theta_\omega\theta_f =I_\mathcal{M}$. Define $ U\coloneqq\theta_g, V\coloneqq\theta_\omega.$ Then $ U:\mathcal{M} \rightarrow\ell^p(\mathbb{N})$ is a Lipschitz  right-inverse of $ \theta_\tau$ and  $V: \ell^p(\mathbb{N}) \rightarrow \mathcal{X}$ is  a linear bounded left-inverse of $ \theta_f$ such that the operator $ VU=\theta_\omega\theta_g=S_{g,\omega}$ is invertible. Further,
	\begin{align*}
		(\zeta_nU)x=\zeta_n\left(\sum_{k=1}^\infty g_k(x)e_k\right)=\sum_{k=1}^\infty g_k(x)\zeta_n(e_k)=g_n(x), \quad \forall x \in \mathcal{M}
	\end{align*} 
	and $Ve_n=\theta_\omega e_n=\omega_n, \forall n \in \mathbb{N} $.
\end{proof}
\begin{lemma}\label{ASFLEMMA2}
	Let $ (\{f_n \}_{n}, \{\tau_n \}_{n}) $ be a  Lipschitz p-ASF for   $\mathcal{M}\subseteq\mathcal{X}$. Then 
	\begin{enumerate}[\upshape(i)]
		\item $R: \mathcal{M} \rightarrow \ell^p(\mathbb{N})$ is a Lipschitz  right-inverse of $ \theta_\tau$  if and only if 
		\begin{align*}
			R=\theta_fS_{f,\tau}^{-1}+(I_{\ell^p(\mathbb{N})}-\theta_fS_{f,\tau}^{-1}\theta_\tau)U
		\end{align*} where $U:\mathcal{M} \to \ell^p(\mathbb{N})$ is a Lipschitz map.
		\item  $ L:\ell^p(\mathbb{N})\rightarrow \mathcal{X}$ is a bounded left-inverse of $ \theta_f$ if and only if 
		\begin{align*}
			L=S_{f,\tau}^{-1}\theta_\tau+V(I_{\ell^p(\mathbb{N})}-\theta_fS_{f,\tau}^{-1}\theta_\tau),
		\end{align*} 
		where $V:\ell^p(\mathbb{N}) \to \mathcal{X}$ is a bounded linear operator. 
	\end{enumerate}		
\end{lemma}
\begin{proof}
	\begin{enumerate}[\upshape(i)]
		\item  $(\Leftarrow)$  $\theta_\tau(\theta_fS_{f,\tau}^{-1}+(I_{\ell^p(\mathbb{N})}-\theta_fS_{f,\tau}^{-1}\theta_\tau)U)=I_\mathcal{M}+\theta_\tau U-I_\mathcal{M}\theta_\tau U=I_\mathcal{M}$. Therefore $\theta_fS_{f,\tau}^{-1}+(I_{\ell^p(\mathbb{N})}-\theta_fS_{f,\tau}^{-1}\theta_\tau)U$ is a Lipschitz  right-inverse of $ \theta_\tau$.  
		
		$(\Rightarrow)$  Define $U\coloneqq R $. Then $\theta_fS_{f,\tau}^{-1}+(I_{\ell^p(\mathbb{N})}-\theta_fS_{f,\tau}^{-1}\theta_\tau)U=\theta_fS_{f,\tau}^{-1}+(I_{\ell^p(\mathbb{N})}-\theta_fS_{f,\tau}^{-1}\theta_\tau)R=\theta_fS_{f,\tau}^{-1}+R-\theta_fS_{f,\tau}^{-1}=R$.
		\item
		$(\Leftarrow)$  $(S_{f,\tau}^{-1}\theta_\tau+V(I_{\ell^p(\mathbb{N})}-\theta_fS_{f,\tau}^{-1}\theta_\tau))\theta_f=I_\mathcal{M}+V\theta_f-V\theta_fI_\mathcal{M}=I_\mathcal{M}$. Therefore  $S_{f,\tau}^{-1}\theta_\tau+V(I_{\ell^p(\mathbb{N})}-\theta_fS_{f,\tau}^{-1}\theta_\tau)$ is a bounded left-inverse of $\theta_f$.
		
		$(\Rightarrow)$  Define $V\coloneqq L$. Then $S_{f,\tau}^{-1}\theta_\tau+V(I_{\ell^p(\mathbb{N})}-\theta_fS_{f,\tau}^{-1}\theta_\tau) =S_{f,\tau}^{-1}\theta_\tau+L(I_{\ell^p(\mathbb{N})}-\theta_fS_{f,\tau}^{-1}\theta_\tau)=S_{f,\tau}^{-1}\theta_\tau+L-S_{f,\tau}^{-1}\theta_\tau= L$.
	\end{enumerate}		
\end{proof}
\begin{theorem}\label{ALLDUAL}
	Let $ (\{f_n \}_{n}, \{\tau_n \}_{n}) $ be a Lipschitz p-ASF for   $\mathcal{M}\subseteq \mathcal{X}$. Then a  Lipschitz p-ASF   $ (\{g_n \}_{n}, \{\omega_n \}_{n}) $ for $\mathcal{M}$  is a dual  for $ (\{f_n \}_{n}, \{\tau_n \}_{n}) $ if and only if
	\begin{align*}
		&g_n=f_nS_{f,\tau}^{-1}+\zeta_nU-f_nS_{f,\tau}^{-1}\theta_\tau U,\\
		&\omega_n=S_{f,\tau}^{-1}\tau_n+Ve_n-V\theta_fS_{f,\tau}^{-1}\tau_n, \quad \forall n \in \mathbb{N}
	\end{align*}
	such that 
	\begin{align*}
		S_{f,\tau}^{-1}+VU-V\theta_fS_{f,\tau}^{-1}\theta_\tau U
	\end{align*}
	is an invertible bi-Lipschitz map, where   $U:\mathcal{M} \to \ell^p(\mathbb{N})$ is a Lipschitz map,  $ V:\ell^p(\mathbb{N})\to \mathcal{X}$ is a  bounded linear operator,  $VU(\mathcal{M})\subseteq \mathcal{M}$, $Ve_n\in \mathcal{M} $ for all $n \in \mathbb{N}$ and 
	\begin{align*}
		&\sum_{n=1}^{\infty}\zeta_n(\theta_fS_{f,\tau}^{-1}+(I_{\ell^p(\mathbb{N})}-\theta_fS_{f,\tau}^{-1}\theta_\tau)U)x)[S_{f,\tau}^{-1}+VU-V\theta_fS_{f,\tau}^{-1}\theta_\tau U]^{-1}(S_{f,\tau}^{-1}\theta_\tau+V(I_{\ell^p(\mathbb{N})}-\theta_fS_{f,\tau}^{-1}\theta_\tau))e_n\\
		&=x\quad \forall x \in \mathcal{M}.
	\end{align*}
\end{theorem}
\begin{proof}
	Lemmas \ref{ASFLEMMA1} and  \ref{ASFLEMMA2} give the characterization of dual frame as 
	\begin{align*}
		&g_n=\zeta_n\theta_fS_{f,\tau}^{-1}+\zeta_nU-\zeta_n\theta_fS_{f,\tau}^{-1}\theta_\tau U=f_nS_{f,\tau}^{-1}+\zeta_nU-f_nS_{f,\tau}^{-1}\theta_\tau U,\\
		&\omega_n=S_{f,\tau}^{-1}\theta_\tau e_n+Ve_n-V\theta_fS_{f,\tau}^{-1}\theta_\tau e_n=S_{f,\tau}^{-1}\tau_n+Ve_n-V\theta_fS_{f,\tau}^{-1}\tau_n, \quad \forall n \in \mathbb{N}
	\end{align*}
	such that 
	$$(S_{f,\tau}^{-1}\theta_\tau+V(I_{\ell^p(\mathbb{N})}-\theta_fS_{f,\tau}^{-1}\theta_\tau))(\theta_fS_{f,\tau}^{-1}+(I_{\ell^p(\mathbb{N})}-\theta_fS_{f,\tau}^{-1}\theta_\tau)U) $$
	is  an invertible bi-Lipschitz map, where $U:\mathcal{M} \to \ell^p(\mathbb{N})$ is a Lipschitz map, $ V:\ell^p(\mathbb{N})\to \mathcal{X}$ is a  bounded linear operator, $VU(\mathcal{M})\subseteq \mathcal{M}$, $Ve_n\in \mathcal{M} $ for all $n \in \mathbb{N}$ and 
		\begin{align*}
		\sum_{n=1}^{\infty}\zeta_n(\theta_fS_{f,\tau}^{-1}+(I_{\ell^p(\mathbb{N})}-\theta_fS_{f,\tau}^{-1}\theta_\tau)U)x)W^{-1}(S_{f,\tau}^{-1}\theta_\tau+V(I_{\ell^p(\mathbb{N})}-\theta_fS_{f,\tau}^{-1}\theta_\tau))e_n=x\quad \forall x \in \mathcal{M}, 
	\end{align*}
where 
\begin{align*}
	W\coloneqq (S_{f,\tau}^{-1}\theta_\tau+V(I_{\ell^p(\mathbb{N})}-\theta_fS_{f,\tau}^{-1}\theta_\tau))(\theta_fS_{f,\tau}^{-1}+(I_{\ell^p(\mathbb{N})}-\theta_fS_{f,\tau}^{-1}\theta_\tau)U).
\end{align*}
	Through expansion and simplification we get 
	\begin{align*}
		(S_{f,\tau}^{-1}\theta_\tau+V(I_{\ell^p(\mathbb{N})}-\theta_fS_{f,\tau}^{-1}\theta_\tau))(\theta_fS_{f,\tau}^{-1}+(I_{\ell^p(\mathbb{N})}-\theta_fS_{f,\tau}^{-1}\theta_\tau)U)
		=S_{f,\tau}^{-1}+VU-V\theta_fS_{f,\tau}^{-1}\theta_\tau U.
	\end{align*}
\end{proof}
 Balan introduced the notion of similarity for frames for Hilbert space which gives an equivalence relation on frames  \cite{BALAN}. It  has been done for Banach spaces by Krishna and Johnson in  \cite{KRISHNAJOHNSON} . We define the same for Lipschitz p-ASF as follows.
  \begin{definition}\label{SIMDEF}
  	Two Lipschitz p-ASFs $ (\{f_n\}_{n}, \{\tau_n\}_{n}) $ and $ (\{g_n \}_{n}, \{\omega_n \}_{n}) $ for $\mathcal{M}\subseteq\mathcal{X}$ are said to be \textbf{similar} or \textbf{equivalent} if there exist invertible bi-Lipschitz map   $T_{f,g}:\mathcal{M} \to \mathcal{M}$ and an invertible bounded linear operator $T_{\tau,\omega} :\mathcal{X} \to \mathcal{X}$ such that $T_{\tau,\omega}(\mathcal{M})\subseteq \mathcal{M}$ and 
  	\begin{align*}
  		g_n=f_nT_{f,g},\quad  \omega_n= T_{\tau,\omega}\tau_n, \quad \forall  n \in \mathbb{N}.
  	\end{align*}
  \end{definition}
  Since  maps giving similarity are invertible,   similarity  is an equivalence relation  on the set 
  \begin{align*}
  	\{(\{f_n\}_{n}, \{\tau_n\}_{n}): (\{f_n\}_{n}, \{\tau_n\}_{n}) \text{ is a  Lipschitz p-ASF for } \mathcal{M}\}.
  \end{align*}
	Observe that for every Lipschitz p-ASF  $(\{f_n\}_{n}, \{\tau_n\}_{n}),$ both  
\begin{align*}
	(\{f_n{S}_{f, \tau}^{-1}\}_{n}, \{\tau_n\}_{n})	 \text{ and }  (\{f_n\}_{n}, \{{S}_{f, \tau}^{-1}\tau_n\}_{n})
\end{align*}   
are  Lipschitz p-ASFs and are  similar to  $(\{f_n\}_{n}, \{\tau_n\}_{n})$.  Balan gave an operator algebraic characterization of similarity in Hilbert spaces \cite{BALAN} and it is extended to Banach spaces by Krishna and Johnson  in \cite{KRISHNAJOHNSON}. We derive Lipschitz version in the following theorem.
  \begin{theorem}\label{SEQUENTIALSIMILARITY}
  	For two Lipschitz p-ASFs $ (\{f_n\}_{n}, \{\tau_n\}_{n}) $ and $ (\{g_n \}_{n}, \{\omega_n \}_{n}) $ for $\mathcal{M}\subseteq\mathcal{X}$, the following are equivalent:
  	\begin{enumerate}[\upshape(i)]
  		\item   $g_n=f_nT_{f, g} , \omega_n=T_{\tau,\omega}\tau_n,  \forall  n \in \mathbb{N}$, for some invertible bi-Lipschitz map $T_{f,g}:\mathcal{M} \to \mathcal{M}$,  for some invertible linear map $T_{\tau,\omega}:\mathcal{X} \to \mathcal{X}$ such that $T_{\tau,\omega}(\mathcal{M})\subseteq \mathcal{M}$. 
  		\item $\theta_g=\theta_f T_{f,g}, \theta_\omega=T_{\tau,\omega}\theta_\tau$,  for some  invertible bi-Lipschitz map $T_{f,g}:\mathcal{M} \to \mathcal{M}$,  for some invertible linear map $T_{\tau,\omega}:\mathcal{X} \to \mathcal{X}$ such that $T_{\tau,\omega}(\mathcal{M})\subseteq \mathcal{M}$. 
  		\item $P_{g,\omega}=P_{f, \tau}.$
  	\end{enumerate}
  	If one of the above conditions is satisfied, then  invertible maps in  $\operatorname{(i)}$ and  $\operatorname{(ii)}$ are unique and given by  $T_{f,g}= S_{f,\tau}^{-1}\theta_\tau\theta_g, T_{\tau, \omega}=\theta_\omega\theta_fS_{f,\tau}^{-1}.$ In the case that $ (\{f_n \}_{n}, \{\tau_n \}_{n}) $ is a Lipschitz p-SF, then $ (\{g_n \}_{n}, \{\omega_n \}_{n}) $ is  a Lipschitz p-SF if and only if $T_{\tau, \omega}T_{f,g} =I_\mathcal{M}$   if and only if $ T_{f,g}T_{\tau, \omega} =I_\mathcal{M}$. 
  \end{theorem}
  \begin{proof}
  	(i) $\Rightarrow $ (ii) 
  	Let $x\in \mathcal{M}$ and $\{a_n\}_{n} \in \ell^p(\mathbb{N})$. Then 
  	\begin{align*}
  	&\theta_gx=\{g_n(x)\}_{n}=\{f_n(T_{f,g}x)\}_{n}=\theta_f(T_{f,g}x),\\
  	&	\theta_\omega(\{a_n\}_{n})=\sum_{n=1}^\infty a_n\omega_n=\sum_{n=1}^\infty a_nT_{\tau,\omega}\tau_n=T_{\tau,\omega}\theta_\tau\{a_n\}_{n}.
  	\end{align*}
  	(ii) $\Rightarrow $ (iii) $  S_{g,\omega}= \theta_\omega\theta_g=T_{\tau,\omega} \theta_\tau\theta_f T_{f,g} =T_{\tau,\omega} S_{f, \tau}T_{f,g}$ and 
  	\begin{align*}
  		P_{g,\omega}=\theta_g S_{g,\omega}^{-1} \theta_\omega=(\theta_f T_{f,g})(T_{\tau,\omega} S_{f, \tau}T_{f,g})^{-1}(T_{\tau,\omega} \theta_\tau)= P_{f, \tau}.
  	\end{align*}  
  	(ii) $\Rightarrow $ (i)  $ \sum_{n=1}^\infty g_n(x)e_n=\theta_g(x)=\theta_f(T_{f,g}x)=\sum_{n=1}^\infty f_n(T_{f,g}x)e_n, \forall x \in \mathcal{M}.$ This  gives (i).\\
  	(iii) $\Rightarrow $ (ii) $\theta_g=P_{g,\omega} \theta_g= P_{f,\tau}\theta_g=\theta_f(S_{f,\tau}^{-1}\theta_{\tau}\theta_g)$ and $\theta_\omega=\theta_\omega P_{g,\omega}=\theta_\omega P_{f,\tau}=(\theta_\omega\theta_fS_{f,\tau}^{-1})\theta_\tau .$ We  show that $S_{f,\tau}^{-1}\theta_{\tau}\theta_g$ and $\theta_\omega\theta_fS_{f,\tau}^{-1} $ are invertible. For,
  	\begin{align*}
  		&(S_{f,\tau}^{-1}\theta_{\tau}\theta_g)(S_{g,\omega}^{-1}\theta_{\omega}\theta_f)=S_{f,\tau}^{-1}\theta_{\tau}P_{g,\omega}\theta_f=S_{f,\tau}^{-1}\theta_{\tau} P_{f,\tau}\theta_f=I_\mathcal{M},\\
  		&(S_{g,\omega}^{-1}\theta_{\omega}\theta_f)(S_{f,\tau}^{-1}\theta_{\tau}\theta_g)=S_{g,\omega}^{-1}\theta_{\omega} P_{f,\tau}\theta_g=S_{g,\omega}^{-1}\theta_{\omega}P_{g,\omega}\theta_g=I_\mathcal{M} 
  	\end{align*}
  	and 
  	\begin{align*}
  		&(\theta_\omega\theta_fS_{f,\tau}^{-1})(\theta_\tau\theta_gS_{g,\omega}^{-1})=\theta_\omega P_{f,\tau}\theta_gS_{g,\omega}^{-1}=\theta_\omega P_{g,\omega}\theta_gS_{g,\omega}^{-1}=I_\mathcal{M},\\
  		&(\theta_\tau\theta_gS_{g,\omega}^{-1})(\theta_\omega\theta_fS_{f,\tau}^{-1})=\theta_\tau P_{g,\omega}\theta_fS_{f,\tau}^{-1}=\theta_\tau P_{f,\tau}\theta_fS_{f,\tau}^{-1}=I_\mathcal{M}.
  	\end{align*}    
  	Let $T_{f,g}, T_{\tau,\omega}:\mathcal{M} \to \mathcal{M}$ be  invertible bi-Lipschitz maps and $g_n=f_nT_{f, g}, \omega_n=T_{\tau,\omega}\tau_n,  \forall  n \in \mathbb{N}$. Then $\theta_g=\theta_fT_{f, g} $ says that $\theta_\tau\theta_g=\theta_\tau\theta_fT_{f, g}=S_{f,\tau}T_{f, g}  $ which implies $ T_{f, g} =S_{f,\tau}^{-1}\theta_\tau\theta_g$. Similarly  $\theta_\omega=T_{\tau,\omega}\theta_\tau $ says that  $\theta_\omega\theta_f=T_{\tau,\omega}\theta_\tau\theta_f=T_{\tau,\omega}S_{f,\tau} $. Hence $T_{\tau,\omega}=\theta_\omega\theta_fS_{f,\tau}^{-1} $. 
  \end{proof}
In Definition \ref{DUALDEF} we defined  the notion of dual frames. There is an  opposite  notion associated is the notion of duality called orthogonality  which is studied for Hilbert spaces in \cite{HANKORNELSONLARSONWEBER, HANMEMOIRS} and for Banach spaces in \cite{KRISHNAJOHNSON}. We can define the orthogonality for Lipschitz p-ASFs as follows.
\begin{definition}\label{ORTHOGONALDEF}
	Let $ (\{f_n\}_{n}, \{\tau_n\}_{n}) $ be a Lipschitz p-ASF for 	$\mathcal{M}\subseteq\mathcal{X}$. 	A Lipschitz p-ASF $ (\{g_n \}_{n}, \{\omega_n \}_{n}) $ for $\mathcal{M}$ is said to be \textbf{orthogonal}   for $ (\{f_n \}_{n}, \{\tau_n \}_{n}) $ if 
	\begin{align*}
		0=\sum_{n=1}^\infty g_n(x) \tau_n=\sum_{n=1}^\infty
		f_n(x) \omega_n, \quad \forall x \in
		\mathcal{M}.
	\end{align*}
\end{definition}
Similar to Proposition \ref{ORTHOGONALPRO} we have the following proposition.
\begin{proposition}
	Given  two Lipschitz p-ASFs $ (\{f_n\}_{n}, \{\tau_n\}_{n}) $ and $ (\{g_n \}_{n}, \{\omega_n \}_{n}) $ for $\mathcal{M}\subseteq\mathcal{X}$, the following are equivalent:
	\begin{enumerate}[\upshape(i)]
		\item  $ (\{g_n \}_{n}, \{\omega_n \}_{n}) $ is  orthogonal  for $ (\{f_n \}_{n}, \{\tau_n \}_{n}) $.
		\item $\theta_\tau\theta_g =\theta_\omega\theta_f =0$.
	\end{enumerate}
\end{proposition}
Using orthogonality  we  derive following  interpolation result. For the Hilbert space frames this is derived by Han and Larson in  \cite{HANMEMOIRS} and for Banach spaces in \cite{KRISHNAJOHNSON}.
\begin{theorem}\label{INTERPOLATION}
	Let $ (\{f_n\}_{n}, \{\tau_n\}_{n}) $ and $ (\{g_n \}_{n}, \{\omega_n \}_{n}) $ be  two Lipschitz p-SF  for  $\mathcal{M}\subseteq\mathcal{X}$ which are  orthogonal. If $A,B, :\mathcal{M}\to \mathcal{M}$ are bi-Lipschitz  maps,   $C, D :\mathcal{X}\to \mathcal{X}$ are bounded linear operators, $C(\mathcal{M})\subseteq \mathcal{M}$, $D(\mathcal{M})\subseteq \mathcal{M}$ and  $ CA+DB=I_\mathcal{M}$, then  
	\begin{align*}
		(\{f_nA+g_nB\}_{n}, \{C\tau_n+D\omega_n\}_{n})
	\end{align*}
	is a  Lipschitz p-SF for  $\mathcal{M}$. In particular,  if scalars $ a,b,c,d$ satisfy $ca+db =1$, then 
	$ (\{af_n+bg_n\}_{n}, \{c\tau_n+d\omega_n\}_{n}) $ is a  Lipschitz p-SF for  $\mathcal{M}$.
\end{theorem} 
\begin{proof}
	We find  
	\begin{align*}
		\theta_{fA+gB} x = \{(f_nA+g_nB)(x) \}_{n}=\{f_n(Ax) \}_{n}+\{g_n(Bx) \}_{n}=\theta_f(Ax)+\theta_g(Bx), \quad \forall x \in \mathcal{M}
	\end{align*}
	and 
	\begin{align*}
		\theta_{C\tau+D\omega}\{a_n \}_{n}=\sum_{n=1}^\infty a_n(C\tau_n+D\omega_n)=C\theta_\tau\{a_n \}_{n}+D\theta_\omega\{a_n \}_{n}, \quad  \forall \{a_n\}_n  \in \ell^p(\mathbb{N}). 
	\end{align*} 
	So 	
	\begin{align*}
		S_{fA+gB,C\tau+D\omega} &=\theta_{C\tau+D\omega} \theta_{fA+gB}= ( C\theta_\tau+ D\theta_\omega)(\theta_fA+\theta_gB)\\
		&=C\theta_\tau\theta_fA+C\theta_\tau\theta_gB+D\theta_\omega\theta_fA+D\theta_\omega\theta_gB\\
		&=CS_{f,\tau}A+0+0+DS_{g,\omega}B
		=CI_\mathcal{M}A+DI_\mathcal{M}B=I_\mathcal{M}.
	\end{align*}
\end{proof} 
We use Theorem \ref{SEQUENTIALSIMILARITY}  to relate three notions duality, similarity and orthogonality.
\begin{proposition}\label{LASTONE}
	For every Lipschitz p-SF $(\{f_n\}_{n}, \{\tau_n\}_{n})$ for $\mathcal{M}\subseteq\mathcal{X}$, the canonical dual for $(\{f_n\}_{n}, \{\tau_n\}_{n})$ is the only dual Lipschitz p-ASF that is similar to $(\{f_n\}_{n}, \{\tau_n\}_{n})$.
\end{proposition}
\begin{proof}
	Let us suppose that  two Lipschitz ASFs $(\{f_n\}_{n}, \{\tau_n\}_{n})$ and $ (\{g_n \}_{n}, \{\omega_n \}_{n}) $  are similar and dual to each other. Then there exist  invertible bi-Lipschitz maps  $T_{f,g}, T_{\tau,\omega} :\mathcal{M}\to \mathcal{M}$  such that $ g_n=f_nT_{f,g},\omega_n=T_{\tau,\omega}\tau_n ,\forall n \in \mathbb{N}$. Theorem \ref{SEQUENTIALSIMILARITY} then gives
	\begin{align*}
		T_{f,g}=S_{f,\tau}^{-1}\theta_\tau\theta_g=S_{f,\tau}^{-1}I_\mathcal{M}=S_{f,\tau}^{-1}\text{ and }T_{\tau, \omega}=\theta_\omega\theta_fS_{f,\tau}^{-1}=I_\mathcal{M}S_{f,\tau}^{-1}=S_{f,\tau}^{-1}.
	\end{align*} 
	Hence $ (\{g_n \}_{n}, \{\omega_n \}_{n}) $ is the canonical dual for  $(\{f_n\}_{n}, \{\tau_n\}_{n})$.	
\end{proof}
\begin{proposition}\label{LASTTWO}
	Two similar Lipschitz p-ASFs cannot be orthogonal.
\end{proposition}
\begin{proof}
	Let $(\{f_n\}_{n}, \{\tau_n\}_{n})$ and $ (\{g_n \}_{n}, \{\omega_n \}_{n}) $  be two Lipschitz p-ASFs for $\mathcal{M}\subseteq\mathcal{X}$ which are similar.  Then there exist   invertible bi-Lipschitz maps  $T_{f,g}, T_{\tau,\omega} :\mathcal{M}\to \mathcal{M}$  such that $ g_n=f_nT_{f,g},\omega_n=T_{\tau,\omega}\tau_n ,\forall n \in \mathbb{N}$.  Theorem \ref{SEQUENTIALSIMILARITY} then says   $\theta_g=\theta_f T_{f,g}, \theta_\omega=T_{\tau,\omega}\theta_\tau  $.    Therefore 
	\begin{align*}
		\theta_\tau \theta_g=\theta_\tau\theta_f T_{f,g}=S_{f,\tau}T_{f,g}\neq0. 
	\end{align*}
\end{proof}
Another use of orthogonal frames is to take direct sum. Given Lipschitz maps $f,g:\mathcal{M} \to \mathbb{K}$, we define $f\oplus g: \mathcal{M} \oplus \mathcal{M} \ni x \oplus y \mapsto f(x)+g(y)
\in \mathbb{K}$.
\begin{theorem}
Let $ (\{f_n\}_{n}, \{\tau_n\}_{n}) $ and $ (\{g_n \}_{n}, \{\omega_n \}_{n}) $ be  two Lipschitz p-ASFs  for  $\mathcal{M}\subseteq\mathcal{X}$ which are  orthogonal.  Then  $
(\{f_n\oplus g_n\}_{n}, \{\tau_n\oplus \omega_n\}_{n})$
is a  Lipschitz p-ASF for  $\mathcal{M}\oplus \mathcal{M}\subseteq \mathcal{X}\oplus \mathcal{X}$. 	
\end{theorem}
\begin{proof}
Let $x\oplus y \in 	\mathcal{M}\oplus \mathcal{M}$. Then 
\begin{align*}
	S_{f\oplus g, \tau\oplus \omega}(x\oplus y )&=\sum_{n=0}^{\infty}(f_n\oplus g_n)(x\oplus y)(\tau_n\oplus \omega_n)\\
	&=\left(\sum_{n=0}^{\infty}f_n(x)\tau_n+\sum_{n=0}^{\infty}g_n(x)\tau_n\right)\oplus \left(\sum_{n=0}^{\infty}f_n(x)\omega_n+\sum_{n=0}^{\infty}g_n(x)\omega_n\right)\\
	&=(S_{f,\tau}x+0)\oplus (0+S_{g, \omega}y)=(S_{f,\tau}\oplus S_{g, \omega})(x\oplus y ).
\end{align*}
\end{proof}

\section{An open problem}
Motivated from the approximation properties of Banach spaces (Schauder basis problem) \cite{ENFLO, CASAZZAAPP} and from the failure of atomic decompositions for (even separable) Banach spaces (see \cite{CASAZZACHRISTENSEN}), we formulate the following interesting (high-end) problem. 
\begin{problem}\label{PROBLEM}
	\begin{enumerate}
		\item  \textbf{Classify which subsets of a Banach space has Lipschitz p-ASF for some $1\leq p<\infty$? In particular, 	does every subset of a Banach space has Lipschitz p-ASF for some $1\leq p<\infty$? }
		\item  \textbf{Classify which subsets $\mathcal{M}$ of a Banach space $\mathcal{X}$  has the following property. There exists a sequence $\{\tau_n\}_n$ in  $\mathcal{M}$ and 	a sequence $\{f_n\}_n$  in  $\operatorname{Lip}(\mathcal{M}, \mathcal{X})$ such that 
 	\begin{align*}
 		x=\sum_{n=1}^\infty f_n(x)\tau_n, \quad \forall x \in
 	\mathcal{M}.	
 	\end{align*}}
 	\end{enumerate}
\end{problem}

\section{Acknowledgments}
 First author thanks Prof. B. V. Rajarama Bhat for the Post-doctoral  position supported by  J. C. Bose
 Fellowship (SERB).

 \bibliographystyle{plain}
 \bibliography{reference.bib}

\end{document}